\documentclass[a4paper,12pt]{article}

\usepackage{makeidx}

\usepackage{latexsym}
\usepackage{amscd}
\usepackage{amsmath}
\usepackage{amssymb}
\usepackage{amsthm}

\usepackage{float}
\usepackage{psfrag}
\usepackage{pst-all}

\usepackage{mymacros}

\psset{nodesep=7pt,labelsep=2pt, arrows=->}

\newcommand{\nflocus}{\ensuremath{X}}

\newcommand{\vecn}{\ensuremath{\vec{n}}}
\newcommand{\vecx}{\ensuremath{\vec{x}}}
\newcommand{\vecy}{\ensuremath{\vec{y}}}
\newcommand{\vecz}{\ensuremath{\vec{z}}}
\newcommand{\vecc}{\ensuremath{\vec{c}}}

\newcommand{\qgr}{\ensuremath{\mathrm{qgr}}}
\newcommand{\bp}{\ensuremath{\boldsymbol{p}}}

\newcommand{\tor}{\ensuremath{\mathrm{tor\text{--}}}}
\newcommand{\gr}{\ensuremath{\mathrm{gr\text{--}}}}
\newcommand{\grproj}{\ensuremath{\mathrm{grproj\text{--}}}}

\newcommand{\dirFukW}{\Fuk^{\rightarrow} W}

\newcommand{\vspan}{\ensuremath{\mathop{\mathrm{span}}}}

\newcommand{\vc}{L}
\newcommand{\Dbsing}{\ensuremath{D^{\mathrm{gr}}_{\mathrm{Sg}}}}

\newcommand{\DG}{\ensuremath{DG}}
\newcommand{\Dg}{\ensuremath{\scD}}
\newcommand{\Tri}{\ensuremath{\scT}}
\newcommand{\PreTr}{\ensuremath{\operatorname{Pre-Tr}}}
\newcommand{\Abp}{A}
\newcommand{\Lbp}{L}

\title{Homological Mirror Symmetry and Simple Elliptic Singularities}
\author{Kazushi Ueda}
\date{}
\begin{document}

\maketitle

\begin{abstract}
We give a full exceptional collection
in the triangulated category of singularities
in the sense of Orlov
\cite{Orlov_DCCSTCS}
for a hypersurface singularity of Fermat type,
and discuss its relation
with homological mirror symmetry
for simple elliptic hypersurface singularities.
\end{abstract}

\section{Introduction}

In \cite{KSaito_EES},
Saito defined
a simple elliptic singularity
to be a normal surface singularity
such that the exceptional set of the minimal resolution
is a smooth elliptic curve,
and classified those
which are hypersurface singularities
into the following three types:
\begin{eqnarray}
 E_6 : x^3 + y^3 + z^3 + \lambda x y z = 0, \label{eq:E6}\\
 E_7 : x^4 + y^4 + z^2 + \lambda x y z = 0, \label{eq:E7}\\
 E_8 : x^6 + y^3 + z^2 + \lambda x y z = 0. \label{eq:E8}
\end{eqnarray}
Here, $\lambda$ is the parameter
which determines the complex structure
of the exceptional curve.
They form a good class of surface singularities
which come next to simple singularities.

Simple singularities are know to be closely related
to simple Lie algebras;
the Milnor lattice of a simple singularity
is isomorphic to the root lattice
of the simple Lie algebra
of the corresponding type,
and the semiuniversal deformation
and its simultaneous resolution
can be constructed
in terms of this Lie algebra
after the works of
Brieskorn, Grothendieck, and Slodowy
(see, e.g., \cite{Brieskorn_SESSAG}).

To generalize this relation
between singularities and Lie algebras
to simple elliptic singularities,
Saito \cite{KSaito_EARSI}
introduced the notion of an {\em elliptic root system}
by abstracting the properties of
vanishing cycles sitting in
the Milnor lattice of a simple elliptic singularity.
The Lie algebra associated with an elliptic root system
was constructed by
Saito and Yoshii \cite{KSaito_EARSIV}.
For
the elliptic root system
coming from the Milnor lattice
of a simple elliptic hypersurface singularity,
this Lie algebra turns out to be
the universal central extension
of $\frakg[s, s^{-1}, t, t^{-1}]$,
where
$\frakg$ is the simple Lie algebra
of the corresponding type.
This universal central extension
was also studied
by Moody, Rao, and Yokonuma
\cite{Moody-Rao-Yokonuma_TLAVR}
under the name of the
{\em $2$-toroidal algebra}.

Recently,
a realization of the ``positive part''
of the quantized enveloping algebra
of an elliptic Lie algebra
in terms of the Ringel--Hall algebra
of the category of coherent sheaves
on a weighted projective line
\cite{Geigle-Lenzing_WPC}
of genus one
was found
by Schiffmann \cite{Schiffmann_NPCQLA}.
The Ringel--Hall algebra
of an abelian category
is an algebra
spanned by isomorphism classes
of indecomposable objects
as a free abelian group,
whose structure constants
are given by ``counting numbers of extensions.''
It was originally used by Ringel \cite{Ringel_HAQG}
to construct the positive part of
the quantized enveloping algebra
of a simply-laced simple Lie algebra
from the abelian category of finite-dimensional
representations of a quiver
of the corresponding type.

There is also a work
of Lin and Peng
\cite{Lin-Peng_ELATA}
who used the Ringel--Hall algebra
of the {\em root category}
of representations of certain quivers
to construct elliptic Lie algebras.
Here, the root category
is the orbit category of the derived category
by twice the shift functor,
which was used by Peng and Xiao
\cite{Peng-Xiao_RCSLA}
to answer the problem
posed by Ringel \cite{Ringel_RARTFDA}
of constructing
not only the positive part
but the whole Lie algebra
from the Ringel-Hall algebra.
The construction of Lin and Peng
is related to that of Schiffmann
in that the derived category of representations
of their quiver is
equivalent
to the derived category of coherent sheaves
on the weighted projective line
used by Schiffmann.

These constructions of elliptic Lie algebras
suggest a link between simple elliptic singularities
and weighted projective lines of genus one.
Elliptic Lie algebras come
from elliptic root systems,
and
in Schiffmann's work,
the elliptic root lattice
is realized
as the Grothendieck group
of the derived category of coherent sheaves
on the weighted projective line.
In this sense,
the derived category of coherent sheaves
on a weighted projective line of genus one
is a {\em categorification}
of an elliptic root lattice.
On the singularity side,
the same elliptic root lattice comes
as the Milnor lattice
of a simple elliptic singularity.
This line of thought naturally leads
to search
for a categorification
of the Milnor lattice.

Such a categorification is provided
by the {\em directed Fukaya category}
defined by Seidel \cite{Seidel_VC}.
It is an $A_\infty$-category
whose set of objects is
a distinguished basis of vanishing cycles
and whose spaces of morphisms are
Lagrangian intersection Floer complexes.
Although the directed Fukaya category
depends on the choice of a distinguished basis
of vanishing cycles,
different choices are related by {\em mutations},
and the derived category
is independent of this choice.

Now a natural question is the relation
between the derived category
of the directed Fukaya category
of a simple elliptic hypersurface singularity
and the derived category of coherent sheaves
on a weighted projective line of genus one.
Since both categories have
full exceptional collections,
we would like
to compare them.
The problem then is that
there are many possible choices
for a full exceptional collection
in a triangulated category.

Let us first discuss the singularity side.
Fix a field $k$.
For a positive integer $p$
greater than one,
let $W_p(X) \in \bC[X]$
be a general polynomial of degree $p$
in one variable.
Then for a suitable choice of
a distinguished basis of vanishing cycles,
the corresponding directed Fukaya category
$\dirFuk W_p$ over $k$
has $p-1$ objects
$( L_i )_{i=1}^{p-1}$,
and the spaces of morphisms
are
$$
 \hom_{\dirFuk W_p}(L_i, L_j) = 
   \begin{cases}
    k \cdot \id_{L_i} & j=i, \\
    k[-1] & j = i + 1, \\
    0 & \text{otherwise},
   \end{cases}
$$
where $k[-1]$ is the one-dimensional graded vector space
concentrated in degree $1$
\cite[section (2B)]{Seidel_VC2}.
Note that 
$\dirFuk W_{p}$ is
not only an $A_\infty$-category
but also a $\DG$-category with a trivial differential,
since all the $A_\infty$-operations
are necessarily trivial.

Now
for positive integers
$p_0, \dots, p_n$
greater than one,
consider the polynomial
\begin{equation} \label{eq:E8W}
 W_{p_0, \ldots, p_n}(X_0, \ldots, X_n)
  = W_{p_0}(X_0) + \cdots + W_{p_n}(X_n)
\end{equation}
in $n+1$ variables,
which defines a holomorphic map
from $\bC^{n+1}$ to $\bC$.
By equipping $\bC^{n+1}$
with the standard symplectic structure,
$W$ gives an exact Lefschetz fibration
in the sense of Seidel \cite{Seidel_LES},
and one can consider its
directed Fukaya category $\dirFuk W_{p_0, \dots, p_n}$.
The following conjecture is a special case
of \cite[Conjecture 1.3]{Auroux-Katzarkov-Orlov_WPP}:
\begin{conjecture}
\label{conj:HMSfP}
There exists an equivalence
\begin{equation} \label{eq:HMSfP}
 D^b \dirFuk W_{p_0, \ldots, p_n}
  \cong D^b(\dirFuk W_{p_0} \otimes \cdots \otimes \dirFuk W_{p_n})
\end{equation}
of triangulated categories.
\end{conjecture}
This conjecture is known to hold
when $n = 1$
\cite[section 6.3]{Auroux-Katzarkov-Orlov_WPP},
and the general case
is under investigation
by Auroux, Katzarkov, Orlov, and Seidel.
With this in mind,
we consider
$D^b(\dirFuk W_{p_0} \otimes \cdots \otimes \dirFuk W_{p_n})$
instead of $D^b \dirFuk W_{p_0, \dots, p_n}$
in this paper,
which clearly has
a full exceptional collection
$
 (L_1 \otimes \dots \otimes L_1, \dots,
   L_{p_0-1} \otimes \dots \otimes L_{p_n-1} ).
$
%

Let us next discuss the weighted projective side.
For a sequence $\bp = (p_0, p_1, p_2)$ of
integers greater than one,
put
\begin{equation} \label{eq:Abp}
 \Abp(\bp)
  = k[x, y, z] / (x^{p_0} + y^{p_1} + z^{p_2}).
\end{equation}
Let $L(\bp)$ be the abelian group
generated by
four elements
$\vecx, \vecy,$, $\vecz$, and $\vecc$
with the relation
$p_0 \vecx = p_1 \vecy = p_2 \vecz = \vecc$,
and equip $\Abp(\bp)$ with an $L(\bp)$-grading
by setting
$\deg x = \vecx$,
$\deg y = \vecy$, and
$\deg z = \vecz$.
\begin{definition}
The category $\qgr \Abp(\bp)$
of coherent sheaves
on the weighted projective line
of weight $\bp = (p_0, p_1, p_2)$
is the quotient category
$$
 \qgr \Abp(\bp) = \gr \Abp(\bp) / \tor \Abp(\bp)
$$
of the category $\gr \Abp(\bp)$
of finitely-generated $L(\bp)$-graded $\Abp(\bp)$-modules
by its full subcategory $\tor \Abp(\bp)$
consisting of torsion modules.
\end{definition}
This definition is equivalent to
that of Geigle and Lenzing
by the Serre's Theorem
in \cite[section 1.8]{Geigle-Lenzing_WPC}.

The first candidate of a full exceptional collection
on a weighted projective line
is the Beilinson-type generator
provided by \cite[Proposition 4.1]{Geigle-Lenzing_WPC}.
However,
this collection
does not match nicely
with the collection
$(L_1 \otimes L_1 \otimes L_1, \dots
  L_{p_0-1} \otimes L_{p_1-1} \otimes L_{p_2-1})$
on the singularity side,
and one needs another
method to find a
full exceptional collection
on the weighted projective line.

Such a method is provided
by a theorem of Orlov \cite{Orlov_DCCSTCS}
relating the derived category of coherent sheaves
to the {\em triangulated category of singularities}.
The origin of this category
goes back to the work of Eisenbud
on matrix factorizations
\cite{Eisenbud_HACI},
which was revived by Kapustin and Li \cite{Kapustin-Li_DBLGM}
and Orlov \cite{Orlov_TCS}
motivated by an idea of Kontsevich.
The concept of grading
was introduced
by Hori and Walcher \cite{Hori-Walcher_FTE},
and the relation
with the derived category of coherent sheaves
was conjectured
by Walcher
\cite[section 4.7]{Walcher_SLGB}.
The idea that
triangulated categories
of singularities
may be useful
in constructing Lie algebras
from hypersurface singularities
is due to Takahashi \cite{Takahashi_MF}
and Kajiura, Saito, and Takahashi
\cite{Kajiura-Saito-Takahashi_MF}.
Although
the above works
deal only with gradings by $\bZ$,
it is straightforward to extend the construction
to gradings by slightly more involved groups
which are extensions of finite abelian groups
by $\bZ$,
and these extra structures
provide a link
with weighted projective lines.

\begin{definition}
The triangulated category of singularities
$\Dbsing(\Abp(\bp))$
is the quotient category
\begin{equation*}
 \Dbsing(\Abp(\bp))
  = D^b(\gr \Abp(\bp)) / D^b (\grproj \Abp(\bp)).
\end{equation*}
of the bounded derived category
$D^b (\gr \Abp(\bp))$
of finitely-generated $L(\bp)$-graded $\Abp(\bp)$-modules
by its full triangulated subcategory
$D^b (\grproj \Abp(\bp))$
consisting of perfect complexes,
i.e.,
bounded complexes of projective modules.
\end{definition}

The following theorem
is obtained
by a straightforward adaptation
of Orlov's argument in \cite{Orlov_DCCSTCS}
to the $L(\bp)$-graded situation:
\begin{theorem} \label{th:DCCSTCS}
For
$\bp=(3,3,3), (2,4,4)$ or $(2,3,6)$,
there exists
an equivalence
$$
 \Dbsing(\Abp(\bp)) \cong D^b(\qgr \Abp(\bp))
$$
of triangulated categories.
\end{theorem}

The advantage of working with $\Dbsing(\Abp(\bp))$
instead of $D^b(\qgr \Abp(\bp))$
is that
one can easily find
an exceptional collection
$(E_1, \ldots, E_{(p_0-1)(p_1-1)(p_2-1)})$
whose total morphism algebra
$
 \bigoplus_{i, j = 1}^{(p_0-1)(p_1-1)(p_2-1)} \Ext^*(E_i, E_j)
$
is isomorphic to
that of
$
 (L_1 \otimes L_1 \otimes L_1, \dots, 
  L_{p_0-1} \otimes L_{p_1-1} \otimes L_{p_2-1}).
$
One can also show that
this exceptional collection is {\em full},
i.e., generates the whole triangulated category,
and that
the $\DG$-algebra underlying this graded algebra
is {\em formal},
i.e.,
quasi-isomorphic to its cohomology.
These facts prove the following:

\begin{theorem} \label{th:main}
For a sequence
$\bp = (p_0, p_1, p_2)$
of integers greater than one,
there exists
an equivalence
$$
 \Dbsing(\Abp(\bp))
  \cong D^b(\dirFuk W_{p_0} \otimes \dirFukW_{p_1} \otimes \dirFuk W_{p_2})
$$
of triangulated categories.
\end{theorem}

By combining Theorem \ref{th:DCCSTCS}
with Theorem \ref{th:main},
one obtains
\begin{equation} \label{eq:HMS}
 D^b(\qgr \Abp(\bp)) \cong
  D^b(\dirFuk W_{p_0} \otimes \dirFuk W_{p_1} \otimes \dirFuk W_{p_2})
\end{equation}
for $(p_0, p_1, p_2) = (3,3,3), (2,4,4)$ and $(2,3,6)$.
The right-hand side of (\ref{eq:HMS}) is expected,
by Conjecture \ref{conj:HMSfP},
to be equivalent to
the derived category $D^b \dirFuk W_{p_0, p_1, p_2}$
of the directed Fukaya category
of $W_{p_0, p_1, p_2}$,
which is a deformation of a simple elliptic singularity.
As a slogan,
{\em the mirrors of the weighted projective lines
of weights $(3,3,3)$, $(2,4,4)$, and $(2,3,6)$
are simple elliptic hypersurface singularities}.
When $\bp = (2,4,4)$ or $(2,3,6)$,
one can state as a theorem
that the mirror of the weighted projective line
of weight $\bp$
is the exact Lefschetz fibration
given by $W_{p_1,p_2}$,
since
Conjecture \ref{conj:HMSfP}
is known to hold for $n=1$
and tensoring with $\dirFuk W_{2}$
does not change the category.
This adds two more examples
to the list
\cite{Auroux-Katzarkov-Orlov_WPP,
Auroux-Katzarkov-Orlov_dP,
Seidel_VC2,
Ueda_HMSTdPS}
of spaces
with a full exceptional collection
where the {\em homological mirror symmetry} conjecture
of Kontsevich
\cite{Kontsevich_HAMS, Kontsevich_ENS98}
is known to hold.

{\bf Acknowledgment}:
We thank
Takeshi Abe,
Takuro Abe,
Yuichiro Hoshi,
and Hiraku Kawanoue
for valuable discussions.
The author is supported by
the 21st Century COE Program of
Osaka University.

{\em Note added} :
We are grateful to
Atsushi Takahashi
who informed us
that he has a more direct proof
of the equivalence
in (\ref{eq:HMS})
based on mutations
of full exceptional collections
along the lines
of Ebeling
\cite{Ebeling_QFMS,Ebeling_MLGBSSS}
and Gabrielov
\cite{Gabrielov_IMCS}.

\section{Generalities on triangulated categories}

In this section,
we recall the definitions of admissible subcategories,
semiorthogonal decompositions,
exceptional collections,
and enhanced triangulated categories
following Bondal and Kapranov \cite{Bondal-Kapranov_ETC}
and Orlov \cite{Orlov_DCCSTCS}.
Throughout this section,
we assume that all categories are small.

\begin{definition}
Let $\scN$
be a full triangulated subcategory
of a triangulated category $\Tri$.
$\scN$ is right (resp. left) admissible
if the inclusion functor
$i : \scN \hookrightarrow \Tri$
has a right (resp. left) adjoint functor.
$\scN$ is admissible
if it is right and left admissible.
\end{definition}

The importance of the concept of admissibility
lies in the following fact:

\begin{lemma}
Let $\scN$ be a full triangulated subcategory
in a triangulated category $\Tri$.
If $\scN$ is right (resp. left) admissible,
then the quotient category
$\Tri / \scN$ is equivalent to $\scN^\bot$
(resp. $\!\ ^\bot \scN$).
\end{lemma}

Here,
$\scN^\bot$ (resp. $\,^\bot \scN$)
denotes the right (resp. left) orthogonal to $\scN$,
i.e.,
the full subcategory
consisting of all objects $M$ such that
$\Hom(N, M) = 0$
(resp. $\Hom(M, N) = 0$)
for all $N \in \scN$.

Admissible subcategories
give {\em semiorthogonal decompositions}
of triangulated categories:

\begin{definition}
A sequence of
full triangulated subcategories
$(\scN_0, \ldots, \scN_n)$
in a triangulated category $\Tri$
is a weak semiorthogonal decomposition
if there exists
a sequence of left admissible subcategories
$\scN_0 = \Tri_0 \subset \Tri_1 \subset \cdots \subset \Tri_n = \Tri$
such that $\scN_p$ is left orthogonal to $\Tri_{p-1}$ in $\Tri_p$.
We write this as $\Tri = \langle \scN_0, \ldots, \scN_n \rangle$.
\end{definition}

A {\em full exceptional collection} is a generator
of a triangulated category
with good properties:

\begin{definition} \label{def:ec}
Let $\Tri$ be a triangulated category
over a field $k$. \\[-5mm]
\begin{enumerate}
 \item An object $E$ of $\Tri$
is exceptional if
$$
 \Ext^i(E, E)
  = \begin{cases}
     k & \text{if $i = 0$}, \\
     0 & \text{otherwise}.
    \end{cases}
$$
 \item An ordered set of objects $(E_0, \ldots, E_n)$
is an exceptional collection if
all the $E_i$'s are exceptional
and
$$
 \Ext^l(E_i, E_j) = 0
$$
for all $i > j$ and $l \in \bZ$.
 \item An exceptional collection $(E_0, \ldots, E_n)$
is full if they generate $\Tri$ as a triangulated category.
\end{enumerate}
\end{definition}

%

Given a full exceptional collection
$(E_0, \dots, E_n)$,
one might hope to reconstruct $\Tri$ from
its total morphism algebra
$$
 \bigoplus_{i,j=0}^n \Ext^*(E_i, E_j).
$$
Unfortunately, this in general is not possible
due to the loss of information
coming from the triangulated structure of $\Tri$,
such as Massey products.
To remedy this situation,
Bondal and Kapranov \cite{Bondal-Kapranov_ETC}
introduced the concept of
an {\em enhancement}
of a triangulated category.
The philosophy of derived categories is
to deal with complexes rather than their cohomologies,
to keep track of the information
which is lost
in passing from complexes to their cohomologies.
To push this philosophy further,
one has to treat
not only objects
but also the spaces of morphisms
as complexes:

\begin{definition}
Let $k$ be a field.
A {\em differential graded category}
over $k$,
or a {\em $DG$-category},
is a category such that
\begin{enumerate}
 \item the set $\hom(E, F)$ of morphisms
  between two objects $E$ and $F$
  is a complex of $k$-vector space, i.e.,
  has a $\bZ$-grading
  $\hom(E, F) = \bigoplus_{k \in \bZ} \hom^k(E, F)$
  and a differential $d$ of degree 1,
 \item for any object $E$,
  the identity morphism $\id_E \in \hom(E, E)$
  satisfies $d(\id_E) = 0$,
  and
 \item the composition
  $$
   \circ : \hom(E, F) \times \hom(F, G) \to \hom(E, G)
  $$
  is a $k$-bilinear map preserving the grading,
  and satisfies the Leibniz rule
  $$
   d(a \circ b) = d(a) \circ b + (-1)^{\deg a} a \circ d(b).
  $$
\end{enumerate}
A {\em $\DG$-functor} between two $\DG$-categories
is a $k$-linear functor
which preserves the gradings
and commutes with the differentials.
\end{definition}

For a nice review on $\DG$-categories,
see e.g., \cite{Keller_ODGC}.
With a \DG-category $\Dg$,
one can associate its cohomology category
$H^0(\Dg)$
whose set of objects is the same as $\Dg$
and whose set $\Hom(E, F)$ of morphisms
is the zeroth cohomology group of $\hom(E, F)$.

Given two $\DG$-categories,
one can consider its tensor product:
\begin{definition}
The tensor product
$\Dg_1 \otimes \Dg_2$
of two $\DG$-categories
$\Dg_1$ and $\Dg_2$
has
the set of of objects
$$
 \Ob(\Dg_1) \times \Ob(\Dg_2)
$$
and the space of morphisms
$$
 \hom_{\Dg_1 \otimes \Dg_2}((E_1, E_2), (F_1, F_2))
  = \hom_{\Dg_1}(E_1, F_1) \otimes \hom_{\Dg_2}(E_2, F_2)
$$
with the composition
$$
 (f \otimes v) \circ (g \otimes w)
  = (-1)^{\deg v \cdot \deg g}(f \circ g) \otimes (v \circ w)
$$
and the differential
$$
 d(f \otimes v) = (d f) \otimes v + (-1)^{\deg f} f \otimes (d v).
$$
\end{definition}

The following definition is due to
Bondal and Kapranov:
\begin{definition}[{\cite[\S1. Definition 1]{Bondal-Kapranov_ETC}}]
A {\em twisted complex} over a \DG-category $\Dg$
is a set
$\{ \{ E_i \}_{i \in \bZ}, \{ q_{ij} \}_{i, j \in \bZ} \}$,
where $E_i$'s are objects of $\Dg$
equal to $0$ for almost all $i$, and 
$q_{ij}$ is an element of $\hom^{i-j+1}(E_i, E_j)$
satisfying
$d q_{ij} + \sum_k q_{kj} \circ q_{ik} = 0$.
A twisted complex
$\{ E_i, q_{ij} \}$
is called {\em one-sided}
if $q_{ij} = 0$ for $i \geq j$.
\end{definition}
We will assume that all twisted complexes
are one-sided henceforth.
Twisted complexes over a \DG-category $\Dg$
form a \DG-category by
$$
 \hom^k(\{ E_i, q_{ij} \}, \{ F_i, r_{ij} \}) =
  \bigoplus_{l+i-i=k} \hom_{\Dg}^l(E_i, F_j)
$$
and
$$
 df = d_{\Dg} f + \sum_m(r_{jm} \circ f + (-1)^{l(i-m+1)} f \circ q_{mi})
$$
for $f \in \hom_{\Dg}^l(E_i, F_j)$.
Let $\Dg^\oplus$ be the $\DG$-category
obtained form $\Dg$ by formally adjoining
finite direct sums,
and $\PreTr(\Dg)$ be the $\DG$-category
of twisted complexes over $\Dg^{\oplus}$.
The cohomology category
$H^0(\PreTr(\Dg))$
will be denoted by $D^b(\Dg)$.
A twisted complex
$
 K = \{ E_i, q_{ij} \} \in \PreTr(\Dg)
$
defines a contravariant $\DG$-functor
from $\Dg$
to the $\DG$-category
of complexes of $k$-vector spaces
by sending $E \in \Ob(\Dg)$
to $\hom_{\PreTr(\Dg)}(E, K)$.
Here, $E$ is considered as a twisted complex
$\{ E , 0 \}$
concentrated at degree zero.
\begin{definition}[{\cite[\S3. Definition 1]{Bondal-Kapranov_ETC}}]
A $\DG$-category $\Dg$ is {\em pretriangulated}
if for every twisted complex
$K \in \Ob(\PreTr(\Dg))$,
the corresponding contravariant $\DG$-functor
is representable.
\end{definition}
By \cite[\S3. Proposition 2]{Bondal-Kapranov_ETC},
the cohomology category $H^0(\Dg)$
of a pretriangulated $\DG$-category $\Dg$
has a natural structure of a triangulated category.

\begin{definition}
An {\em enhancement}
of a triangulated category $\Tri$
is a pretriangulated $\DG$-category $\Dg$
together with an equivalence
$H^0(\Dg) \to \Tri$
of triangualted categories.
\end{definition}
A typical example of an enhancement
is the $\DG$-category
underlying the bounded derived category
of an abelian category
with enough injectives
or projectives.
\cite[\S3. Example 3]{Bondal-Kapranov_ETC}.

Now we can state the reconstruction of
a triangulated category
from a finite number of generators:

\begin{theorem}
[{Bondal--Kapranov \cite[\S4. Theorem 1]{Bondal-Kapranov_ETC}}]
\label{th:Bondal-Kapranov}
Assume that a triangulated category $\Tri$
generated by a finite number of objects
$E_0, \ldots, E_n$
has an enhancement $\Dg$.
Let $\scA$ be the full $\DG$-subcategory of $\Dg$
consisting of $E_0, \ldots, E_n$.
Then one has an equivalence
$$
 \scT \cong D^b(\scA)
$$
of triangulated categories.
\end{theorem}

\section{Triangulated categories of singularities}
\label{sc:tcs}

We prove Theorem \ref{th:DCCSTCS}
in this section.
The proof is a straightforward adaptation
of arguments from \cite{Orlov_DCCSTCS}
to the $L(\bp)$-graded situation,
which we include for the reader's convenience.
For an $\Lbp(\bp)$-graded
$\Abp(\bp)$-module $M$ and $\vecn \in L(\bp)$,
$M(\vecn)$ will denote the graded $\Abp(\bp)$-module
obtained by shifting the grading by $\vecn$;
$
 M(\vecn)_{\vec{m}}=M_{\vec{n}+\vec{m}}.
$
Put $\bp = (p_0, p_1, p_2) = (3,3,3), (2,4,4)$, or $(2,3,6)$,
$\Abp = \Abp(\bp)$,
and $L = L(\bp)$.
Define
$
 \phi: L \to \bZ
$
by
$$
 \phi(\vecx) = \frac{p_2}{p_0}, \ 
 \phi(\vecy) = \frac{p_2}{p_1}, \ 
 \phi(\vecz) = 1,
$$
and
let $\scS_{<0}$ and $\scS_{\geq 0}$
be the full triangulated subcategories
of $D^b(\gr \Abp)$
generated by
$k(\vec{n})$
for $\phi(\vec{n}) > 0$
and $\phi(\vec{n}) \leq 0$
respectively.
Here $k$ denotes the $\Abp$-module
$\Abp / (x, y, z)$.
Let further
$\scP_{<0}$ and $\scP_{\geq 0}$
be the full triangulated subcategories
of $D^b(\gr \Abp)$
generated by the modules
$\Abp(\vec{n})$
for $\phi(\vec{n}) > 0$
and $\phi(\vec{n}) \leq 0$
respectively.
Then one has weak semiorthogonal decompositions
\begin{eqnarray}
 D^b(\gr \Abp) &= \langle \scS_{<0}, D^b(\gr \Abp_{\geq 0}) \rangle, \\
 D^b(\gr \Abp) &= \langle D^b(\gr \Abp_{\geq 0}), \scP_{<0} \rangle,
\end{eqnarray}
where $\gr \Abp_{\geq 0}$ is the full subcategory
of $\gr \Abp$
consisting of $L(\bp)$-graded $\Abp$-modules $M$
such that
$M_{\vec n} = 0$
if $\phi(\vec n) < 0$.
The corresponding statement
in the $\bZ$-graded case is
\cite[Lemma 2.3]{Orlov_DCCSTCS},
whose proof also works here verbatim.
This shows that
$$
 \ ^\bot \scS_{<0} = \scP_{<0}^\bot.
$$
Since $A$ is Gorenstein,
the duality functor
$$
 D =\bR \Hom_\Abp(-, \Abp):
  D^b(\gr \Abp)^{\circ} \to D^b(\gr \Abp)
$$
is an equivalence
and maps
$
 (\scP_{\leq 0})^{\circ}
$
to
$
 \scP_{\geq 0}.
$
Moreover,
since
$$
\bR \Hom_\Abp(k,\Abp)
 = k(\vecx+\vecy+\vecz-\vecc)[-2]
$$
and
$$
 \phi(\vecx+\vecy+\vecz-\vecc) = 0
$$
in our case,
$D$ also gives an equivalence between
$
 (\scS_{\leq 0})^{\circ}
$
and
$
 \scS_{\geq 0}.
$
Therefore,
$\scS_{\geq 0}$ is right admissible in $D^b(\gr \Abp)$,
$\scP_{\geq 0}$ is left admissible in $D^b(\gr \Abp)$,
and
$$
 \scS_{\geq 0}^\bot = \ ^\bot \scP_{\geq 0}.
$$
Moreover,
$\scS_{<0}$
and
$\scP_{\geq 0}$
are mutually orthogonal.
Hence
one has
weak semiorthogonal decompositions
$$
 D^b(\gr \Abp)
  = \langle \scS_{<0}, \scD_0, \scS_{\geq 0} \rangle 
$$
and
$$
 D^b(\gr \Abp)
  = \langle \scS_{<0}, \scP_{\geq 0}, \scT_0 \rangle,
$$
such that
$
 \scD_0 \cong D^b(\qgr \Abp)$
and
$
 \scT_0 \cong \Dbsing(\Abp).
$
Now one has
\begin{eqnarray*}
 D^b(\gr \Abp)
  &= \langle \scS_{<0}, \scD_0, \scS_{\geq 0} \rangle \\
  &= \langle \scP_{\geq 0}, \scS_{<0}, \scD_0 \rangle \\
  &= \langle \scS_{<0}, \scP_{\geq 0}, \scD_0 \rangle \\
  &= \langle \scS_{<0}, \scP_{\geq 0}, \scT_0 \rangle,
\end{eqnarray*}
which shows
$
 \scT_0 = \scD_0
$
and hence
$$
 \Dbsing(\Abp) \cong D^b(\qgr \Abp).
$$
\qed

\section{A full exceptional collection}

We prove Theorem \ref{th:main}
in this section.
Fix any weight $\bp=(p_0, p_1, p_2)$
and put $\Abp = \Abp(\bp)$
and $\Lbp = L(\bp)$.
We will find a full triangulated subcategory
$\scT$ of $D^b(\gr \Abp)$
equivalent to $\Dbsing(\Abp)$
such that
$\left( k(\vec{n}) \right)_{\vec{n} \in I}$
is a full exceptional collection
in $\scT$,
where
$k = \Abp / (x, y, z)$ and
$$
I = \{a \vecx + b \vecy + c \vecz \in \Lbp \suchthat 
 - p_0 + 2 \leq a \leq 0,
 - p_1 + 2 \leq b \leq 0,
 - p_2 + 2 \leq c \leq 0
\}.
$$
Let $\Lbp_{+}$ be the subset of $\Lbp$ defined by
$$
 \Lbp_{+}
  = \{ - 2 \vec{c} + a \vec{x} + b \vec{y} + c \vec{z} |
        a \geq 1, b \geq 1, c \geq 1 \},
$$
and $\Lbp_{-}$ be the complement $\Lbp \setminus \Lbp_+$.
Let further
$\scS_{-}$ and $\scP_{+}$
be the full triangulated subcategories of $D^b(\gr \Abp)$
generated by
$k(\vec{n})$ for $\vec{n} \in L_{+}$ and 
$\Abp(\vec{m})$ for $\vec{m} \in L_{-}$
respectively.
Then $\scS_{-}$ and $\scP_{+}$ are left admissible
in $D^b(\gr \Abp)$,
and since $\scP_{+} \subset \ ^\bot \scS_{-}$,
one has a weak semiorthogonal decomposition
$$
 D^b(\gr \Abp)
  = \langle \scS_{-}, \scP_{+}, \scT \rangle
$$
such that $\scT \cong \Dbsing(\Abp)$.
One can see that
$k(\vec{n})$
for $\vec{n} \in I$
belongs to $\scT$,
since
$$
 \bR \Hom(k(\vec{m}), k(\vec{n})) = 0
$$
if
$
 \vec{m} \notin \vec{n} + \bN \vec{x} + \bN \vec{y} + \bN \vec{z},
$
and
$$
 \bR \Hom(k(\vec{m}), \Abp(\vec{n})) = 0
$$
if
$
 \vec{m} \neq - \vec{c} + \vec{x} + \vec{y} + \vec{z} + \vec{n}.
$
The $\bR \Hom$'s between them
can be calculated by the following free resolution\\
\begin{align*}
\begin{CD}
 \cdots
  @>>>
 \begin{array}{c}
   \Abp(-\vec{c} - \vec{y} - \vec{z}) \\
    \oplus \\
   \Abp(-\vec{c} - \vec{x} - \vec{z}) \\
    \oplus \\
   \Abp(-\vec{c} - \vec{x} - \vec{y}) \\
    \oplus \\
   \Abp(-2 \vec{c}) \\
 \end{array}
  @>d_4>>
 \begin{array}{c}
   \Abp(-\vec{x} - \vec{y} - \vec{z}) \\
    \oplus \\
   \Abp(-\vec{c} - \vec{z}) \\
    \oplus \\
   \Abp(-\vec{c} - \vec{y}) \\
    \oplus \\
   \Abp(-\vec{c} - \vec{x}) \\
 \end{array}
  @>d_3>>
 \begin{array}{c}
   \Abp(- \vec{y} - \vec{z}) \\
    \oplus \\
   \Abp(-\vec{x} - \vec{z}) \\
    \oplus \\
   \Abp(-\vec{x} - \vec{y}) \\
    \oplus \\
   \Abp(-\vec{c}) \\
 \end{array}
\end{CD} \\
\begin{CD}
  @>d_2>>
 \begin{array}{c}
   \Abp(- \vec{z}) \\
    \oplus \\
   \Abp(- \vec{y}) \\
    \oplus \\
   \Abp(- \vec{x}) \\
 \end{array}
  @>d_1>>
 \Abp
  @>>>
 k,
\end{CD}
\end{align*}
where
\begin{align}
 d_1 &=
  \begin{pmatrix}
   z & y & x
  \end{pmatrix}, \\
 d_2 &=
  \begin{pmatrix}
   -y & -x & 0 & z^{p_2-1} \\
   z & 0 & -x & y^{p_1-1} \\
   0 & z & y & x^{p_0-1}
  \end{pmatrix}, \\
 d_3 &=
  \begin{pmatrix}
   x & -y^{p_1-1} & z^{p_2-1} & 0 \\
   -y & -x^{p_0-1} & 0 & z^{p_2-1} \\
   z & 0 & -x^{p_0-1} & y^{p_1-1} \\
   0 & z & y & x
  \end{pmatrix}, \\
 d_4 &=
  \begin{pmatrix}
   x^{p_0-1} & -y^{p_1-1} & z^{p_2-1} & 0 \\
   -y & -x & 0 & z^{p_2-1} \\
   z & 0 & -x & y^{p_1-1} \\
   0 & z & y & x^{p_0-1}
  \end{pmatrix},
\end{align}
to be
\begin{align*}
 \bR &\Hom(k(a_1 \vec{x} + b_1 \vec{y} + c_1 \vec{z}),
      k(a_2 \vec{x} + b_2 \vec{y} + c_2 \vec{z})) \\
  & =
  \begin{cases}
   k & (a_2, b_2, c_2) = (a_1, b_1, c_1), \\
   k[-1] & (a_2, b_2, c_2)
               = (a_1-1, b_1, c_1),
                 (a_1, b_1-1, c_1),
                 (a_1, b_1, c_1-1), \\
   k[-2] & (a_2, b_2, c_2)
                      = (a_1-1, b_1-1, c_1),
                        (a_1-1, b_1, c_1-1),
                        (a_1, b_1-1, c_1-1),\\
   k[-3] & (a_2, b_2, c_2)
                      = (a_1-1, b_1-1, c_1-1), \\
     0   & \text{otherwise}.
  \end{cases}
\end{align*}
Here, $k[i]$ for $i \in \bZ$
denotes the one-dimensional vector space
concentrated in degree $-i$.
Therefore
$
 \left( k(\vec{n})\right)_{\vec{n} \in I}
$
is an exceptional collection.
It is straightforward to read off
the structure of the Yoneda products
from the above resolution
to show that
the the full subcategory of $\scT$
consisting of $(k(\vecn))_{\vecn \in I}$
is isomorphic
as a graded category to
$
 \dirFuk W_{p_0} \otimes \dirFuk W_{p_1} \otimes \dirFuk W_{p_2}.
$
Moreover,
$\scT$
has an enhancement
induced from that of $D^b(\gr \Abp)$,
and one can see using the above resolution again
that the total morphism $\DG$-algebra
$
 \bigoplus_{\vec{m}, \vec{n} \in I}
  \hom(k(\vec{m}), k(\vec{n})) 
$
of the above exceptional collection
is {\em formal}, i.e., quasi-isomorphic to its cohomology
with the trivial differential
as a $\DG$-algebra.
This shows that
$
 D^b(\dirFuk W_{p_0} \otimes \dirFuk W_{p_1} \otimes \dirFuk W_{p_2})
$
is equivalent to the full triangulated subcategory of $\Dbsing(A(\bp))$
generated by the above exceptional collection.

Now we prove that the image of
$
 (k(\vec{n}))_{\vec{n} \in I}
$
in $\Dbsing(\Abp)$ is full.
Let
$
 G_0 = \Spec k[\Lbp / \bZ \vecc]
$
be the subgroup of
$
 G = \Spec k[\Lbp]
$
isomorphic to
$
 (\bZ / p_0 \bZ) \times (\bZ / p_1 \bZ) \times (\bZ / p_2 \bZ).
$
Since $\Lbp$ is the group of characters of $G$,
the $\Lbp$-grading of $\Abp$ defines
an action of $G$ on $\Abp$,
and
an $\Lbp$-graded $\Abp$-module
is the same thing as a $G$-equivariant $\Abp$-module.
For a $G_0$-module $M$,
let $M^{G_0}$ denote its $G_0$-invariant part.
The restriction
of the action of $G$ on $\Spec \Abp$
to $G_0$
is given by
$$
\begin{array}{cccc}
 G_0 \ni ([i], [j], [k]) : & \Spec \Abp & \longrightarrow & \Spec \Abp \\
 &  \rotatebox{90}{$\in$} & & \rotatebox{90}{$\in$} \\
 & (x, y, z) & \longmapsto
  & (\zeta_{p_0}^i x, \zeta_{p_1}^j y, \zeta_{p_2}^k z),
\end{array}
$$
where $\zeta_n = \exp[2 \pi \sqrt{-1} / n]$
for a positive integer $n$.
Let
$\nflocus \subset \Spec \Abp$
be the closed subscheme
defined by the ideal $(x y z)$
where the action of $G$ is not free.

By \cite[Lemma 1.11]{Orlov_TCS},
an element of $\Dbsing(\Abp)$ is isomorphic
to the shift $M[i]$
of a finitely-generated $\Lbp$-graded $\Abp$-module $M$
by some integer $i$.
Let $M_{x y z}$ be the localization of $M$
by $x y z \in \Abp$.
Then $M_{x y z}^{G_0}$
is an $\Abp_{x y z}^{G_0}$-module
and one has
$
 M_{x y z}^{G_0} \otimes_{\Abp_{x y z}^{G_0}} \Abp_{x y z}
   \cong M_{x y z},
$
since the action of $G_0$
on $\Spec A_{x y z}$ is free.
Take a finitely-generated $\Abp^{G_0}$-module $N$
such that
$
 N_{x^{p_0} y^{p_1} z^{p_2}}
  \cong (M_{x y z})^{G_0}.
$
Then one has an isomorphism
$
 (N \otimes_{\Abp^{G_0}} \Abp)_{x y z}
  \to M_{x y z},
$
which can be extended to a morphism
$
 N \otimes_{\Abp^{G_0}} \Abp \to M
$
by multiplying a power of $x y z$
if necessary.
Since
$
 \Abp^{G_0} = k[x^{p_0}, y^{p_1}, z^{p_2}] / (x^{p_0} + y^{p_1} + z^{p_2})
$
is regular and
$\Abp$ is flat over $\Abp^{G_0}$,
$N \otimes_{\Abp^{G_0}} \Abp$ is perfect.
Therefore, by replacing $M$
with the kernel and the cokernel of the above morphism,
one can assume that $M$ is supported
on $\nflocus$.

Since $\nflocus$ is the union of
$G$-invariant subschemes
$\nflocus_0$,
$\nflocus_1$, and
$\nflocus_2$
defined by $x=0$,
$y=0$, and
$z=0$ respectively,
one can assume that
the support of
$M$ is contained in the subscheme
$\nflocus_0 \cong \Spec \Abp / (x)$.
(In general,
a coherent sheaf supported on the union $S_1 \cup S_2$
of closed subschemes
can be obtained from coherent sheaves
supported on $S_1$ and $S_2$
by taking cones.)
For a finitely-generated $\Abp$-module $M$
supported on $\nflocus_0$,
let $l(M)$ be the minimal integer $n$
such that $x^n M = 0$.
Then one has an exact sequence
$$
 0 \to M' \to M \to M / x M \to 0
$$
such that $M / x M$ is an $\Abp/(x)$-module
and $l(M') < l(M)$.
By replacing $M$ with $M'$
and continuing this process,
one can assume that
$M$ is an $\Abp/(x)$-module.

Since
$x$ is not a zero divisor of $\Abp$,
$\Abp/(x)$ is perfect by the following free resolution:
$$
 0 \to \Abp \to \Abp \to \Abp/(x) \to 0.
$$
Therefore,
a perfect $\Abp/(x)$-module
is also perfect as an $\Abp$-module.
By replacing $\Abp$
with $\Abp/(x) = k[y,z] / (y^{p_1} + z^{p_2})$
and repeating the same argument,
one can see that
any $\Lbp$-graded $\Abp / (x)$-module
can be obtained from $\Abp / (x, y)$-modules
by taking cones up to perfect complexes.
Since the same is true for $\Abp / (y)$-modules
and $\Abp / (z)$-modules,
$\Dbsing(\Abp)$
is generated
by modules supported at $x=y=z=0$
as a triangulated category.

Since any module supported at $x=y=z=0$
can be obtained from
$
 k = \Abp / (x, y, z)
$
by extensions,
it is enough to show that
all sheaves
of the form $k(\vec{m})$ for $\vec{m} \in \Lbp$
can be obtained from
$(k(\vecn))_{\vecn \in I}$
by taking cones up to perfect complexes.
To do this,
first note that the exact sequences
\begin{eqnarray*}
  0 \to k(-\vecz) \to & k[z] / (z^2) & \to k \to 0, \\
  0 \to k(-2 \vecz) \to & k[z] / (z^3) & \to k[z] / (z^2) \to 0, \\
  & \vdots & \\
  0 \to k(-(p_2 - 1) \vecz) \to & k[z] / (z^{p_2}) &
   \to k[z] / (z^{p_2 - 1}) \to 0 
\end{eqnarray*}
of $\Abp$-modules
show that
$k(-(p_2 - 1) \vecz)$ can be obtained
from $k, k(-\vecz), \dots, k(-(p_2-2)\vecz)$
by taking cones up to the perfect module
$k[z]/(z^{p_2}) \cong \Abp / (x, y)$.
Then by shifting the degrees,
one can see that
for any $\vecn \in \Lbp$,
$k(\vecn)$ can be obtained from
either
$
 k(\vecn - \vecz), k(\vecn - 2 \vecz), \dots, k(\vecn - (p_2 - 1) \vecz)
$
or
$
 k(\vecn + \vecz), k(\vecn + 2 \vecz), \dots, k(\vecn + (p_2 - 1) \vecz)
$
by taking cones up to perfect complexes.
The same is true for $\vecx$ and $\vecy$,
which proves that $(k(\vecn))_{\vecn \in I}$ is full.

%

\bibliographystyle{plain}
\bibliography{bibs}

Department of Mathematics,
Osaka University,
Machikaneyama 1-1,
Toyonaka,
Osaka,
560-0043,
Japan.

{\em e-mail address}\ : \  kazushi@cr.math.sci.osaka-u.ac.jp

\end{document}